\documentclass{amsart}
\usepackage{amsmath,amssymb,amsfonts,enumerate,amsthm, amscd}
\newtheorem{thm}{Theorem}[section]

\newtheorem{lem}[thm]{Lemma}
\newtheorem{prop}[thm]{Proposition}

\newtheorem{conj}[thm]{Conjecture}
\begin{document}
\bibliographystyle{amsplain}
%\date{}
\title[groups of prime power order]{Commutativity pattern of finite non-abelian $p$-groups determine their orders}
\author{A. Abdollahi}
\address{Department
 of Mathematics, University of Isfahan, Isfahan 81746-73441,
Iran; School of Mathematics, Institute for Research in Fundamental
Sciences (IPM), P.O.Box: 19395-5746, Tehran, Iran.}
\email{a.abdollahi@math.ui.ac.ir}
\author{S. Akbari}
\address{Department of Mathematical Sciences, Sharif University of Technology, Tehran, Iran; School of Mathematics, Institute for Research in Fundamental Sciences (IPM), P.O.Box: 19395-5746, Tehran, Iran.}
\email{s\_akbari@sharif.edu}
\author{H. Dorbidi}
\address{Department of Mathematical Sciences, Sharif University of Technology, Tehran, Iran.}
\email{hr\_dorbidi@yahoo.com}
\author{H. Shahverdi}
\address{Department  of Mathematics, University of Isfahan, Isfahan 81746-73441,
Iran.} \email{hamidshahverdi@gmail.com}
\thanks{The first author's research was in part supported by a grant from IPM (No. 90050219) as well as by the Center of Excellence for Mathematics, University of Isfahan. The second author's research was in part supported by a grant from IPM (No. 90050212)}
\begin{abstract}
Let $G$ be a non-abelian group and $Z(G)$ be the center of $G$.
Associate a graph $\Gamma_G$ (called non-commuting graph of $G$)
with $G$ as follows: take $G\setminus Z(G)$ as the vertices of
$\Gamma_G$ and join two distinct vertices $x$ and $y$, whenever
$xy\neq yx$. Here, we prove that ``the commutativity pattern of a finite non-abelian $p$-group determine its order among the class of groups"; this means that
 if $P$ is a finite
non-abelian $p$-group such that $\Gamma_P\cong \Gamma_H$ for some
group $H$, then $|P|=|H|$.
\end{abstract}
\maketitle
\section{\bf Introduction and Results}
Given a finite non-abelian group $G$, one can associate in many
different ways a graph to $G$ (e.g. \cite{BHM,W}).   Here we
consider the non-commuting graph $\Gamma_G$ of $G$:  the set of
vertices of $\Gamma_G$ is $G\setminus Z(G)$, and two vertices $x$
and $y$ are adjacent if and only if $xy\neq yx$. The non-commuting
graph was first considered by Paul Erd\"os in 1975 \cite{N}. The
non-commuting graph of finite groups has been studied by many
people (e.g., \cite{AAM,M}).\\
 The non-commuting graph of a group is a discrete way to reflect the commutativity pattern of the group.
In \cite{AAM} the following conjecture was  formulated:
\begin{conj}[Conjecture 1.1 of \cite{AAM}]\label{con1}
Let $G$ and $H$ be two finite non-abelian groups such that
$\Gamma_G\cong \Gamma_H$. Then $|G|=|H|$.
\end{conj}
 Conjecture \ref{con1} was refuted in \cite{M} by exhibiting two groups $G$ and $H$ of orders
$$|G| = 2^{10}\cdot 5^3 \not= 2^3 \cdot 5^6 = |H|$$
with isomorphic non-commuting graphs.

In \cite{AAM}, it is proved that Conjecture \ref{con1} holds
whenever one of the groups in question is a symmetric group,
dihedral group, alternative group or a non-solvable AC-group
(where by an AC-group we mean a group in which the centralizer of
every non-central element is abelian). Recently Darafsheh \cite{D}
has proved the validity of Conjecture \ref{con1} whenever one of
the groups $G$ or $H$ is a non-abelian finite simple group.

The main result of the present paper shows that any pair of groups consisting a counterexample for Conjecture \ref{con1} cannot contain a group of prime power order.
\begin{thm}\label{mainthm} If $P$ is a finite
non-abelian $p$-group such that $\Gamma_P\cong \Gamma_G$ for some
group $G$, then $|P|=|G|$.
\end{thm}
This is a curious general phenomenon for  non-abelian groups of prime power order: the order of a  prime power order group can be determined among all finite groups by a proper model of its commutativity behavior, i.e, the non-commuting graph.
\section{\bf{Preliminary Results}}
 It is not hard to prove that the finiteness or the being non-abelian of a group can be transferred
 under graph isomorphism whenever two groups have the same non-commuting graph.
 Throughout  $P$ denotes a fixed but arbitrary finite non-abelian $p$-group of order $p^n$ whose center $Z(P)$ is
 of order $p^r$ and  $1<p^{a_1}<p^{a_2}<\cdots<p^{a_k}$ are all distinct
 conjugacy class sizes of $P$, where $p^{a_i}$ is the size of conjugacy  class ${g_i}^G$ of the element $g_i$.
 Throughout we also denote by $u$ the greatest common divisor $\gcd(a_1,\dots,a_k,n-r)$ of $\{a_1,\dots,a_k,n-r\}$.
\begin{lem}\label{l1} Let $G$ be a finite non-abelian group and $H$ be a group such that
$\phi:\Gamma_G\rightarrow\Gamma_H$ is a graph isomorphism. Then
the following hold:
\begin{enumerate}
\item  $|C_H(h)|$ divides $(|g^G| - 1) (|Z(H)| - |Z(G)|)$, where
$h=\phi(g)$.
\item If $|Z(G)|\geq|Z(H)|$ and  $G$ contains a
non-central element  $g$ such that ${|C_G(g)|}^2\geq
|G|\cdot|Z(G)|$, then $|G|=|H|$.
\end{enumerate}
\end{lem}
\begin{proof}
(1) \; Since $\Gamma_G\cong \Gamma_H$, we have
  $$ |G| - |Z(G)| =|H|-|Z(H)| \Rightarrow |H| = |G|-|Z(G)|+|Z(H)| \eqno{(a)}$$ and
  $$ |C_G(g)|-|Z(G)|=|C_H(h)|-|Z(H)|  \Rightarrow  |C_H(h)|=|C_G(g)|+|Z(H)|-|Z(G)| \eqno{(b)}$$

As $|C_H(h)|$ divides $|H|$, $|C_H(h)|$ divides
$$ (|C_G(g)|+|Z(H)|-|Z(G)|)\frac{|G|}{|C_G(g)|}, \eqno{(c)}$$ it follows from $(a),(b),(c)$ $|C_H(h)|$
 divides $$(|{g^G}|-1)(|Z(H)|-|Z(G)|).$$
(2) \; Let $h=\phi(g)$. By part (1), we have
$|C_H(h)|=|C_G(g)|+|Z(H)|-|Z(G)|$ divides
$(|{g^G}|-1)(|Z(H)|-|Z(G)|)$. Now, the inequality
${|C_G(g)|}^2\geq |G||Z(G)|$ implies that
$$0\leq |C_H(h)|\leq (|g^G|-1)(|Z(G)|-|Z(H)|)<|C_G(g)|+|Z(H)|-|Z(G)|=|C_H(h)|$$
and
this yields  $(|{g^G}|-1)(|Z(G)|-|Z(H)|)=0$. Hence
$|Z(G)|=|Z(H)|$.
\end{proof}
\begin{lem}\label{l2} Suppose that $H=P_1\times A$ is a finite group,  where $P_1$ is a $p$-group, $A$ is a finite abelian
 group such that $\gcd(p,|A|)=1$. If $\Gamma_P\cong \Gamma_H$, then $|P|=|H|$.
\end{lem}
\begin{proof}
  Let $\phi$ be a graph isomorphism from $\Gamma_P$ to $\Gamma_H$.
Suppose  $h=\phi(g_t)$ for some $1\leq t\leq k$ and $|P_1|=p^{\kappa}, |Z(P_1)|=p^{\omega}, |A|=a$ and $|C_H(h)|=ap^{\nu}$.
Since $\Gamma_P\cong\Gamma_H$, we have
$$|P|-|Z(P)|=p^r(p^{n-r}-1)=ap^\omega(p^{\kappa -\omega}-1)=|H|-|Z(H)|,$$ $$|P|-|C_P(g_t)|=p^{n-a_t} (p^{a_t}-1)=ap^{\nu}(p^{\kappa - \nu}-1)=|H|-|C_H(h)|,$$
since $\gcd(a,p)=1$, it follows that $r=\omega$ and $ n-a _t = \nu$. Therefore $$|C_P(g_t)|-|Z(P)|=p^r(p^{n-a_t-r}-1)=ap^{\omega}(p^{\nu-\omega}-1)=|C_H(h)|-|Z(H)|.$$ Therefore
$a=1$. Since $r=\omega$, $|Z(P)|= |Z(H)|$. Hence $|P|=|H|$.
\end{proof}
\begin{lem} \label{thmnew} Suppose $H=Q\times A$, where $Q$ is a $q$-group for some prime $q$,  $A$ is an abelian group and $\gcd(|A|,q)=1$. If $\Gamma_P\cong \Gamma_H$,  then $|H|=|P|$.
\end{lem}
\begin{proof} If $p=q$, then Lemma \ref{l2} completes the proof.
Suppose, for a contradiction, that  $p\neq q$.

 Note that   $|g_1^P|=p^{a_1}$. Let $\phi$ be a graph isomorphism from $\Gamma_P$ to $\Gamma_H$ and   let
$$\phi(g_1)=h,|A|=a ,|Q|=q^\kappa, |C_H
(h)|=aq^{\nu},|Z(H)|=aq^\omega. $$ It is clear that $\kappa
>\nu>\omega$. Since $\Gamma_P\cong\Gamma_H$, we
have
\begin{equation} |C_H(h)|-|Z(H)|=aq^{\omega}(q^{\nu-\omega}-1)=p^r(p^{n-a_1-r}-1)=|C_P(g_1)|-|Z(P)|,\end{equation}
\begin{equation} |H|-|C_H(h)|=aq^{\nu}(q^{\kappa-\nu}-1)=p^{n-a_1}(p^{a_1}-1)=|P|-|C_P(g_1)|.\end{equation}
Since $|g_1^P|\leq |g^P|$ for all $g\in P\setminus Z(P)$,
 $h^H$ has the minimum size  among all conjugacy classes of non-central elements of $H$. By considering the conjugacy class equation of $H$, we have
 $$aq^{\kappa}=aq^{\omega}+q^{\kappa-\nu}+\sum_{i=1}^s|{x_i}^H|,$$
 where $$\{{x_i}^{H}\mid i=1,\dots, s\}=\{g^H \;|\; g\in H\setminus Z(H)\}\setminus\{ h^H\}.$$  Since $\gcd(a,q)=1$ and $q^{\kappa-\nu}|\sum_{i=1}^s|{x_i}^H|$, it follows that
$$\kappa-\nu\leq\omega. \eqno{(*)}$$
Equation (1) implies that the largest  $p$-power number possibly dividing  $a$ is $p^r$.
Now it follows from  Equations (1), (2) and the inequality $(*)$ that
$$ p^{n-a_1-r}|q^{\kappa-\nu}-1\leq q^{\omega}-1\leq
p^{n-a_1-r}-2, $$ which is a contradiction. This completes the proof.
\end{proof}
\begin{lem}\label{t1} Let $H$  be a group such that $\Gamma_P\cong \Gamma_H$. Then
 $|Z(H)|$ divides $p^r(p^u-1)$, where $u=\gcd(a_1,\dots,a_k,n-r)$.
\end{lem}
\begin{proof}
(1) \; Since $\Gamma_P\cong\Gamma_H$,  $|P|-|Z(P)|=|H|-|Z(H)|$ and
$|P|-|C_P(g_i)|=|H|-|C_H(h_i)|$, for every $i\in \{{1,\dots,k}\}$ and $h_i=\phi(g_1)$, where $\phi:\Gamma_G\rightarrow \Gamma_H$.
Therefore we have the following equalities
$$p^r(p^{n-r}-1)=|Z(H)|\big(\frac{|H|}{|Z(H)|}-1\big)$$
$$p^{n-a_i}(p^{a_i}-1)=|C_H(h_i)|\big(\frac{|H|}{|C_H(h_i|)}-1\big)$$ for each $i\in\{1,\dots,k\}$. Thus $|Z(H)|$ divides the great common divisors
of the left hand side of two latter equalities which is $p^r(p^u-1)$.
\end{proof}
A class of groups  arising  in the proof of our main result is the class of AC-groups; as we mentioned, a group $G$ is called an AC-group whenever the centralizer of every non-central element is abelian.
AC-groups was studied by many people (e.g., \cite{S}).
It is easy to see that $C_G(x) \cap C_G(y)=Z(G)$ for any two non-central elements $x,y\in G$ with distinct centralizers.
This implies that $$\mathfrak{C}(G)=\big\{C_G(x)/Z(G) \;|\; x\in G\setminus Z(G)\big\}$$ is a
{\em partition} of         $G/Z(G)$; where by a partition  for a group $H$ we mean a collection $\mathcal{C}$ of
proper subgroups of $H$ such that $H=\bigcup_{S\in \mathcal{C}} S$
and $S \cap T=1$ for any two distinct $S,T\in \mathcal{C}$. Each
element of $\mathcal{C}$  is called a component of the partition.
If each component is abelian, we call  $\mathcal{C}$ an  abelian
partition. Thus $\mathfrak{C}(G)$  is an abelian partition for ${G}/{Z(G)}$.
The size of $\mathfrak{C}(G)$ is an invariant of the non-commuting graph $\Gamma_G$, called the clique number; where by definition
the clique number of a finite graph is the maximum number of
vertices which are pairwise adjacent. The clique number of the
non-commuting graph $\Gamma_H$ of a non-abelian group $H$ will be
denoted by $\omega(H)$. Thus $\omega(H)$ is simply  the
maximum number of pairwise non-commuting elements in the group.
\begin{lem}\label{l6} Suppose that  $G$ is a finite non-abelian AC-group such that ${G}/{Z(G)}$ is a  $p$-group. Then $\omega(G)\equiv 1 \mod p$.
\end{lem}
\begin{proof} Since $G$ is an AC-group, $\omega=\omega(G)=|\mathfrak{C}(G)|$, where $$\mathfrak{C}(G)=\big\{C_G(x) \;|\; x\in G\setminus Z(G)\big\}.$$
On the other hand, $C_G(x) \cap C_G(y)=Z(G)$ for any two
non-central elements $x,y\in G$ such that $C_G(x)\not=C_G(y)$.
Therefore $$|G|=-(\omega-1)|Z(G)|+\sum_{S\in \mathfrak{C}(G)} |S|.$$ This completes the proof.
\end{proof}
\begin{lem}\label{mann}[Mann \cite{ber}: Lemma 39.8, p. 354] Suppose $C$ is a subgroup of group $G$ and let $a\in G$ be such that $CC^a=C^aC$. Then $CC^a=C[C,a]$.
\end{lem}
\begin{proof} We have $$CC^a=\bigcup_{c\in C}Cc^a=\bigcup_{c\in C}Cc^{-1}c^a=\bigcup_{c\in C}C[c,a]\subseteq C[C,a].$$
Thus $CC^a\subseteq C[C,a]$. Since all the generators  $[c,a]=c^{-1}c^a$ ($c\in C$) of $[C,a]$ belongs to $CC^a$, we have $C[C,a]\subseteq CC^a$, since by hypothesis  $CC^a$ is a group. This completes the proof.
\end{proof}
In the following proposition we will use this property of any AC-groups $G$; for any two commuting non-central elements $x$ and $y$ of $G$, we have $C_G(x)=C_G(y)$.
\begin{prop} \label{bigt} Let $G$ be  a nilpotent AC-group of nilpotency class greater than $2$, then the set $\mathfrak{C}$ of all centralizers of non-central elements of  $G$ has exactly one normal member $T$ in $G$. In particular, $T$ is a characteristic subgroup of $G$.   Moreover, the latter normal subgroup $T$ has the maximum order among all members of $\mathfrak{C}$.
\end{prop}
\begin{proof}
Let $x$ be any element of $Z_2(G)\setminus Z(G)$. Then $C_G(x)$ is a normal subgroup
of $G$ containing $G'$: for the map  $\phi$ defined
on $G$ by $g^\phi=[x,g]$ for all $g\in G$ is a group homomorphism
and its image is contained in $Z(G)$ and its kernel is $C_G(x)$.
Since $G$ is of nilpotency class greater than $2$, there exists an element $g\in G'\setminus Z(G)$. Since $[Z_2(G),G']=1$, the remark preceding the proposition implies that $$C_G(x)=C_G(g)  \;\;\text{for all} \;\; x\in Z_2(G)\setminus Z(G). \eqno{\diamondsuit}$$
Now suppose that  $N=C_G(y)$ is  a normal
centralizer of $G$ for some non-central element $y$. Then there exists an element
 $t\in \big(N \cap Z_2(G)\big)\setminus Z(G)$, since $Z(G)\lneqq N$. Since $yt=ty$, it follows from $\diamondsuit$ that $C_G(t)=C_G(y)=C_G(g)$.
Hence, we have so far proved that $\mathfrak{C}$ has exactly one normal member in $G$. This implies that $C_G(x)$ is a characteristic subgroup of $G$.

Now, we prove $C_G(x)$ has the maximum order among all members of $\mathfrak{C}$. Suppose
that $C=C_G(h)$ for some $h\in G\setminus Z(G)$. We may assume that $C$ is not normal in $G$. Thus there exists  an element $a\in
N_G(N_G(C))\setminus N_G(C)$, since $G$ is nilpotent.  Then $C^a\neq C$, and $C^a$ is a subgroup of $N_G(C)$.
Let $A=CC^a$.  By Lemma \ref{mann}, we have $$CC_G(x)\supseteq
CG'Z(G)\supseteq C[C,a]Z(G)=CC^aZ(G)=CC^a=A.$$ It follows that
$$\frac{|C||C_G(x)|}{|Z(G)|}=|CC_G(x)|\geq|A|=|CC^a|=\frac{{|C|}^2}{|Z(G)|}$$
 Thus $|C_G(x)|\geq |C|$. This completes the proof.
\end{proof}
The proof of existence of  unique normal centralizer is due to Rocke
\cite[Lemma 3.8]{R}; the argument to prove the existence of a normal centralizer of
maximal order is due to Mann \cite[Theorem 39.7, p. 354]{ber}.
He has proved among all abelian subgroups of maximal order in a
metabelian $p$-group, there exists a normal subgroup. The latter was first proved by  Gillam \cite{Gi}.
\begin{lem}\label{b} Let $P$ be of nilpotency class $2$. Then $a_i\leq r$ for every $i$.
\end{lem}
\begin{proof}
Since $P$ is of nilpotency class  $2$, for every $x\in P\setminus Z(P)$
with class size $p^{a_i}$, the conjugacy class of $x$ is contained
in $xP'\subseteq xZ(P)$. Hence  $ p^{a_i}\leq p^r$. This completes the proof.
\end{proof}
 Now we  will need the following two well known results about Frobenius groups.
\begin{prop}
\begin{enumerate}
\item (see e.g., Theorem 6.7 of \cite{I}) Let $N$ be a normal subgroup of a finite group $G$, and
suppose that $C_G(n)\subseteq N$ for every non-identity element
$n\in N$. Then $N$ is complemented in $G$, and if $1<N<G$, then
$G$ is a Frobenius group with kernel $N$.
\item  (see e.g., Lemma 6.1 of \cite{I}) Let $H$ be a Frobenius group with the kernel $F$ and a complement
$K$, then $|K|$ divides $|F|-1$.
\end{enumerate}
\end{prop}
\begin{lem} \label{e} Let $H=KF$ be a Frobenius group with the kernel
$F$ and a complement $K$. Suppose $1\subset F_1\subseteq F$ is a
normal subgroup of $H$. Then $H_1=KF_1$ is a Frobenius group with
the kernel $F_1$ and a complement $K$.
\end{lem}
\begin{proof} For every non-identity element $f_1$ of $F_1$,  $$C_{H_1}(f_1)=C_H(f_1)\bigcap H_1\subseteq F\bigcap H_1=F\bigcap F_1K=F_1,$$
by the Dedekind modular law. Therefore $H_1$ is a
Frobenius group with the kernel $F_1$. It is clear  that $K$ is a
complement for $F_1$ in $H_1$.
\end{proof}

\section{\bf{Proof of the Main Result}}
In this section we prove our main result, Theorem \ref{mainthm}.

%Suppose that $P$ is a finite non-abelian $p$-group of order $p^n$ whose center $Z(P)$ is
 %of order $p^r$ and  $1<p^{a_1}<p^{a_2}<\cdots<p^{a_k}$ are all distinct
% conjugacy class sizes of $P$, where $p^{a_i}$ is the size of conjugacy  class ${g_i}^G$ of the element $g_i$.
% We also denote by $u$ the greatest common divisor $\gcd(a_1,\dots,a_k,n-r)$ of $\{a_1,\dots,a_k,n-r\}$.
%Suppose also that $G$ is a group such that $\Gamma_P\cong \Gamma_G$. Let $\phi$ be a graph isomorphism from $\Gamma_P$ to $\Gamma_G$.\\

We argue by induction on the order of $P$. If $|P|=p^3$, then $|P|=|G|$ by Proposition 3.20 of \cite{AAM}.
If $P$ is not an AC-group, there exists a non-central element $x\in P$ such that $C_P(x)$ is non-abelian.
If $y=\phi(x)$, then $\Gamma_{C_P(x)}\cong \Gamma_{C_G(y)}$. Now induction hypothesis implies that $|C_P(x)|=|C_G(y)|$ and since
$|P|-|C_P(x)|=|G|-|C_G(y)|$, we have $|P|=|G|$. Thus, we may assume that $P$ is an AC-group and so $G$ is also an AC-group.
By Proposition 3.14 of \cite{AAM}, we may assume that $G$ is solvable.
Therefore by the classification of non-abelian solvable AC-groups in \cite{S}, $G$ is isomorphic to one  the following groups $H_i$ ($i=1,\dots,5$):
\begin{enumerate}
\item $H_1$ is non-nilpotent and it has an abelian normal subgroup
$N$ of prime index and $\omega(H_1)=|N:Z(H_1)|+1$.
\item $H_2/Z(H_2)$ is a
Frobenius group with the Frobenius kernel and complement $F/Z(H_2)$
and $K/Z(H_2)$, respectively and $F$ and $K$ are abelian subgroups
of $H_2$ and $\omega(H_2)=|F:Z(H_2)|+1$.
\item $H_3/Z(H_3)\cong S_4$ and $V$
is a non-abelian subgroup of $H_3$ such that $V/Z(H_3)$ is the Klein
4-group of $H_3/Z(H_3)$ and $\omega(H_3)=13$, where $S_4$ is the
symmetric group of on $4$ letters.
\item $H_4=A\times Q$, where $A$
is an abelian subgroup and $Q$ is an $AC$-group of prime power
order.
\item $H_5/Z(H_5)$ is a Frobenius group with the Frobenius
kernel and complement $F/Z(H_5)$ and $K/Z(H_5)$, respectively and $K$
is an abelian subgroup of $H$. $Z(F)=Z(H_5)$, and $F/Z(H_5)$ is of
prime power order and $\omega(H_5)=|F:Z(H_5)|+\omega(F)$.
\end{enumerate}
By Lemmas 3.11 and 3.12 of \cite{AAM} and Lemma \ref{thmnew}, we may assume that $G$ is isomorphic to either $H_1$  or $H_5$.
Suppose that $G\cong H_1$. Then, obviously $\Gamma_P\cong \Gamma_{H_1}$.
Since $N$ is abelian,  there exists  $h\in H_1\setminus Z(H_1)$
such that $C_{H_1}(h)=N$. As $P$ is an AC-$p$-group, it follows from
Lemma \ref{l6} that   $$\omega(P) \equiv 1 \mod p.$$ Since
$\Gamma_P\cong\Gamma_{H_1}$, we have $$\omega(H_1)=|C_{H_1}(h):Z(H_1)|+1
\equiv 1 \mod p,$$ and so $p\mid |C_{H_1}(h):Z(H_1)|$.  On the other
hand Lemma \ref{l1}(1) implies that, $|C_{H_1}(h)|$ divides
$(p^{a_t}-1)(p^r-|Z(H_1)|)$, where $g_t$ maps to $h$ under a graph
isomorphism from $\Gamma_P$ to $\Gamma_{H_1}$. Thus $p$ divides $|Z(H_1)|$ and so  $p^2 \mid
|C_{H_1}(h)|$. This follows that $p^2$ divides $|Z(H_1)|$. By continuing
this latter process, one obtains that  $p^r$ divides $|Z(H_1)|$ and
so $|Z(H_1)|\geq |Z(P)|$. Now, let $y\in H_1\setminus C_{H_1}(h)$ so that
$H_1=C_{H_1}(h)C_{H_1}(y)$ and
$$|H_1||Z(H_1)|=|C_{H_1}(h)||C_{H_1}(y)|\leq \max\{|C_{H_1}(h)|^2,|C_{H_1}(y)|^2\}.$$ Now, Lemma \ref{l1}(2) implies that $|P|=|H_1|$.

Thus, it remains to deal with the case $G\cong H_5$. Let $H=H_5$
and note that $\Gamma_P\cong \Gamma_H$. We need to introduce some
new notation for the group $H$. Since ${F}/{Z(F)}$ is a
$q$-group for some prime $q$, we set $|F|=bq^{\kappa}$, for some
positive integer $b$ such that $\gcd(b,q)=1$ and therefore
$|Z(H)|=bq^{\omega}$ and $|C_F(f_i)|=|C_H(f_i)|=bq^{\nu_i}$ for
some $f_i\in F\setminus Z(F)$.(Recall that $Z(F)=Z(H)$ in this case) Since
$F$ is nilpotent and non-abelian, we have
$1\leq\omega<\nu_i<\kappa$. Since
$\gcd(|{K}/{Z(H)}|,|{F}/{Z(H)}|)=1$, we have
$|C_H(h)|=|K|=aq^{\omega}$ for some $h\in H\setminus F$ and for
some positive integer $a$. It is clear that $b\mid a$ and
$\gcd(a,q)=1$. Therefore $|H|=aq^{\kappa}$. Suppose that under a graph isomorphism from $\Gamma_H$
to $\Gamma_P$, $h$ maps to $g_t$ for
some integer $1\leq t\leq k$ and $f_i$ maps to $g_i$, where $1\leq
i\leq k$ and $i\neq t$.  Here note that $f_t$ is not defined. Suppose further that  $\beta=a_t$.

We need to prove the following (a), (b), (c) and (d).

{\bf (a)} \; $p\neq q$.\\

{\bf (b)} \; if $p^l$ divides $a$, for some integer $l$, then $p^l$
divides $b$ and  $p^{r+1}$ does not divide $a$. This simply means that the largest $p$-power part of $a$ and $b$ are the same and $p^r$
is the  largest  $p$-power possibly  dividing $a$.\\

{\bf (c)} \; $\Gamma_F$ is a regular graph so that  there exists integers $\nu$ and $\alpha$ such that
$\nu_i=\nu$ and $a_i=\alpha$ for  all $1\leq i\leq k$ and $i\neq t$.\\

{\bf (d)} \; $\nu\leq 2\omega$ and $\kappa\leq 3\omega$.\\

{\bf Proof of (a)} \; Suppose $p=q$. Since $\Gamma_P\cong \Gamma_H$,
$ap^{\omega}(p^{\kappa-\omega}-1)=p^{n-\beta}(p^{\beta}-1)$ and
$bp^{\omega}(p^{\nu_i-\omega}-1)=p^r(p^{n-a_i-r}-1)$. Therefore
$n-\beta=\omega=r$, a contradiction.

{\bf Proof of (b)} \; Since $\Gamma_P\cong \Gamma_H$, we have
\begin{equation}(a-b)q^{\omega}=p^r(p^{n-\beta-r}-1).\end{equation}
Thus $p^r\mid a-b$. This proves the first part of   {\bf (b)} for all
$l\in\{1,\dots,r\}$. Now, suppose $t>r$ and $p^t$ divides $a$ and $p^t\nmid b$.
Equation (3) shows  $p^{r+1}\nmid b$. Now let $i\in\{1,\dots, k\}$ such that $i\not=t$. Then by  the
graph isomorphism, we have
$$p^{n-a_t}-p^{n-a_i}=aq^{\omega}-bq^{\nu_i}. \eqno{(**)}$$
 Since $r+1\geq n-a_t$ and $r+1\geq n-a_i$, it follows from $(**)$ that  $p^{r+1}$ divides $b$, a contradiction.
Now Equation (3) implies that $p^{r+1}\nmid a$ and since $b$ divides $a$, the proof of part (b) follows.

{\bf Proof of (c)} \; Suppose  $\Gamma_F$ is not regular.
Therefore $F$ has two
 centralizers $C_H(f_{i_1})$ and $C_H(f_{i_2})$ of order $bq^{\nu_{i_1}}$ and $bq^{\nu_{i_2}}$, respectively, where $\nu_{i_1}\neq \nu_{i_2}$.
 We may assume that   the  conjugacy class of $f_{i_1}$ in $F$ is of  minimum size among all
 conjugacy classes of non-central elements of $F$. We distinguish
 two cases to reach a contradiction.

(I) Suppose that $\nu_{i_1}-\nu_{i_2}\leq \omega$.
\begin{equation} p^{n-a_{i_2}}-p^r =bq^{\nu_{i_2}}-bq^{\omega}\end{equation}
\begin{equation} p^{n-a_{i_1}}-p^{n-a_{i_2}}=bq^{\nu_{i_1}}-bq^{\nu_{i_2}}\end{equation}
Now it follows from Equations (4),(5) and part {\bf (b)}  that
$$p^{n-a_{i_2}-r}|q^{\nu_{i_1}-\nu_{i_2}}-1\leq q^{\omega}-1\leq p^{n-a_{i_2}-r}-2,$$
a contradiction.

(II) Suppose that $\nu_{i_1}-\nu_{i_2}> \omega$.

We claim that the nilpotency class of $F$ is greater than $2$. If not, then
 Lemma \ref{b}  implies that $$\kappa-\nu_{i_2}\leq
\omega. \eqno{\clubsuit}$$ Since $\nu_{i_1}-\nu_{i_2}>\omega$, $\clubsuit$ is a
contradiction. Therefore the nilpotency class of $F$ is greater than $2$. Since $H$ is an AC-group, $F$ is also an $AC$-group. Therefore every maximal abelian subgroup of $F$ is
centralizer of non-central element of $F$. By Proposition \ref{bigt},
$F$ has a characteristic  centralizer $C_F(f_j)$ of order
$bq^{\nu_j}=bq^{\nu_{i_1}}$ having the maximum order among the proper centralizers. Thus $\nu_{i_1}=\nu_j$ and so
$a_{i_1}=a_j$.  Since $F$ is normal subgroup of $H$, $C_F(f_j)$ is
normal in $H$. Since ${H}/{Z(H)}$ is Frobenius group, by Lemma
\ref{e} ${K}/{Z(H)}{C_F(f_j)}/{Z(H)}$ is a Frobenius group
with the kernel ${C_F(f_j)}/{Z(H)}$ and a complement
${K}/{Z(H)}$. Thus
$$\frac{a}{b}|q^{\nu_{i_1}-\omega}-1.\eqno{\heartsuit}$$ By the graph isomorphism,
we have
\begin{equation}
bq^{\omega}(q^{{\nu_{i_1}}-\omega}-1)=p^r(p^{n-a_{i_1}-r}-1),\end{equation}
\begin{equation}bq^{{\nu_{i_1}}}(\frac{a}{b}q^{\kappa-\nu_{i_1}}-1)=p^{n-a_{i_1}}(p^{a_{i_1}}-1).\end{equation}
Since $\gcd(\frac{a}{b},p)=1$, Equations $\heartsuit$ and (6) imply that
$\frac{a}{b}q^{\omega}|p^{n-a_{i_1}-r}-1$. Equation (7) imply
that $p^{n-a_{i_1}-r}|\frac{a}{b}q^{\kappa-\nu_{i_1}}-1$ and by
the conjugacy class equation $\kappa-\nu_{i_1}\leq \omega$. Therefore
$\frac{a}{b}q^{\omega}<\frac{a}{b}q^{\omega}$, a contradiction.

{\bf Proof of (d)} \; Since $\Gamma_F$ is regular and $F$ is an AC-group, we have
$\omega(F)=\frac{q^{\kappa-\omega}-1}{q^{\nu-\omega}-1}$.
Therefore $\nu-\omega$ divides $\kappa-\omega$. Now, by considering the conjugacy class equation of $F$, we find that
$\nu-\omega\leq \omega$ and $\kappa\leq 3\omega$.\\

 Now we have  two different possibilities on the centralizer orders of $H$:\\

(I) \; $bq^{\nu}>aq^{\omega}$.
 Since $\Gamma_P \cong \Gamma_H$, we have $$p^{n-\alpha}-p^{n-{\beta}}=bq^{\nu}-aq^{\omega},$$
 where  $\beta=a_t$. It follows from the latter equation, Lemma \ref{t1}  and parts {\bf (b)},{\bf (d)} that
$$p^{n-\beta-r}|q^{\nu-\omega}-\frac{a}{b}<q^{\omega}|p^{u}-1< p^{n-\beta-r},$$
a contradiction.

(II) \; $aq^{\omega}>bq^{\nu}$.

 Since $\Gamma_P \cong \Gamma_H$, we have
\begin{equation}aq^{\omega}-bq^{\nu}=p^{n-\beta}-p^{n-\alpha}.\end{equation}
We consider two cases: \\

(i) \; $u<n-\alpha-r$. Since $u\mid n-\alpha-r$,  $2u\leq
n-\alpha-r$.  Since $H/Z(H)$ is a Frobenius group, $|K/Z(H)|=a/b$ divides $|F/Z(F)|-1$.
 Now it follows from  parts {\bf (b)} and {\bf (d)},  Lemma \ref{t1}(1) and Equation (8), we have
$$p^{n-\alpha-r}\mid \frac{a}{b}-q^{\nu-\omega}\leq
q^{\kappa-\omega}-1-q^{\nu-\omega}<q^{2\omega}|(p^u-1)^2<p^{2u},$$
a contradiction.\\

 (ii) \;  $u=n-\alpha-r$. Since $u\mid n-\beta-r$,  $n-\beta-r\geq
2u$. By the graph isomorphism
$$p^n-p^{n-\beta}=aq^{\omega}(q^{\kappa-\omega}-1).$$
This latter
equation, Lemma \ref{t1} and parts {\bf (b)} and {\bf (d)} imply that
$$p^{n-\beta-r}\mid q^{\kappa-\omega}-1<q^{2\omega}\mid (p^u-1)^2<p^{2u},$$
 a contradiction.

This completes the proof. $\hfill \Box$


\begin{thebibliography}{99}
\bibitem{AAM}
A.~Abdollahi, S.~Akbari, H.R.~Maimani, \emph{Non-commuting graph
of a group}, J. Algebra \textbf{298} (2006) 468-492.
\bibitem{ber}
Y.~Berkovich, \emph{Groups of prime power order}, Vol. 1,  de
Gruyter Expositions in Mathematics, 46, Walter de Gruyter GmbH \& Co. KG, Berlin, 2008.
\bibitem{BHM}
B.A.~Bertram, M.~Herzog, A.~Mann, \emph{On a graph related to
conjugacy classes of groups}, Bull. London Math. Soc. {\bf 22} (6)
(1990) 569-575.
\bibitem{Gi} J.D.~Gillam, \emph{A note on finite metabelian p-groups}, Proc.
Amer. Math. Soc. {\bf 25} (1970), 189–190.
\bibitem{D} M.R. Darafsheh, \emph{Groups with the same non-commuting graph},
Discrete Appl. Math. \textbf{157} (2009) no. 4, 833-837.
\bibitem{I}
I.M.~Isaacs, \emph{Finite group theory}, Graduate Studies in Mathematics, 92, Amer. Math. Soc., Providence, RI, 2008.
\bibitem{M}
A.R.~Moghaddamfar,  \emph{About Noncommuting graphs}, Siberian
Math. J. \textbf{47} (2005)  no.5, 1112-1116.
\bibitem{N}
B.H.~Neumann, \emph{A problem of Paul Erd\"os on groups}, J. Aust.
Math. Soc. Ser. A21 (1976) 467-472.
\bibitem{R}D.M.~Rocke,  \emph{$p$-groups with abelian centralizers}, Proc. London Math. Soc. (3) {\bf 30} (1975) 55-75.
\bibitem{S}
R.~Schmidt, \emph{Zenralisatorverb\"ande endlicher Gruppen}, Rend.
Sem. Mat. Univ. Padova \textbf{44} (1975) 55-75.
\bibitem{W}
J.S.~Williams, \emph{Prime graph components of finite groups}, J.
Algebra {\bf 69}  (1981) 487–513.
\end{thebibliography}
\end{document}